\documentclass{amsart}

\usepackage{graphicx, psfrag}

\newtheorem{theorem}{Theorem}
\newtheorem{lemma}[theorem]{Lemma}

\begin{document}

\title{a characterization of hyperbolic geometry among Hilbert geometry}

\author{Ren Guo}

\address{Department of Mathematics, Rutgers University, Piscataway, NJ, 08854}

\email{renguo@math.rutgers.edu}

\thanks{This work is partially supported by NSF Grant \#0604352.}



\begin{abstract}
In this paper we characterize hyperbolic geometry among Hilbert geometry by the property that
three medians of any hyperbolic triangle all pass through one point.
\end{abstract}

\maketitle

\begin{figure}[htbp]
\begin{center}
\includegraphics[scale=.4]{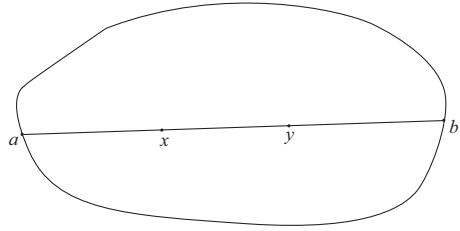}
\end{center}
\caption{\label{hilbert}Hilbert distance}
\end{figure}

We begin by recalling the Hilbert geometry of an open convex set.
Let $K$ be a nonempty bounded open convex set in $\mathbf{R}^n,
n\geq2.$ The \it Hilbert distance \rm $d_K$ on $K$ is introduced
by David Hilbert as follows. For any $x\in K,$ let $d_K(x,x)=0;$
for distinct points $x,y$ in $K$, assume the line passing through
$x,y$ intersects the boundary $\partial \overline{K}$ at two
points $a,b$ such that the order of these four points on the line
is $a,x,y,b$ as in Figure \ref{hilbert}.

Denote the cross-ratio of the points by
$$(a,x,y,b)=\frac{|y-a|}{|x-a|}\cdot\frac{|x-b|}{|y-b|}$$ where
$|\cdot|$ is the Euclidean norm of $\mathbf{R}^n.$ Then the
Hilbert distance is $$d_K(x,y)=\frac12ln(a,x,y,b).$$ The metric
space $(K,d_K)$ is called a Hilbert geometry. Notice that a
straight line in $(K,d_K)$ is a geodesic. When $K$ is the unit
open ball $\{(x_1,...x_n)\in \mathbf{R}^n|\sum_{i=1}^nx^2_i<1\}$,
$(K,d_K)$ is the Klein model of the hyperbolic geometry.

Since the cross-ratio is invariant under any projective mapping
$P$, $(K,d_K)$ and $(P(K),d_{P(K)})$ are isometric as Hilbert
geometry. In particular, when $K$ is an open ellipsoid (ellipse
when n=2), $\{(x_1,...x_n)\in
\mathbf{R}^n|\sum_{i=1}^n\frac{x^2_i}{a^2_i}<1\}$ for some nonzero
numbers $a_1,...a_n$, $(K,d_K)$ is isometric to the hyperbolic
geometry.

A natural question is to characterize hyperbolic geometry among
Hilbert geometry. By the argument above, it is equivalent to
characterize ellipsoid among open convex set. In \cite{cvv}, the
authors characterize ellipsoid by using the fact that every
hyperbolic ideal triangle has a constant area. In this paper we
will consider another elementary property. The following is a
well-known fact in constant sectional curvature geometry (see
\cite{ge} Chapter 7, especially problem K-19):

Three medians of any triangle all pass through one point.

That this property characterizes hyperbolic geometry among Hilbert
geometry was conjectured by Feng Luo \cite{l}. Our main result is
the affirmative solution of it.

\begin{theorem}\label{1} Let $K$ be a nonempty bounded open convex set in $\mathbf{R}^n$ with a Hilbert distance
$d_K$. If three medians of any triangle in $(K,d_k)$ all pass through one point (called property M), then $K$ is an open ellipsoid.
\end{theorem}

\begin{proof} First we consider the case $n=2.$
If $K$ is not an open ellipse, we try to find a triangle which
dose not have the property M. We will use the following lemma
which is a simplified version of a result due to Fritz John
\cite{j}, see also \cite{ba} \cite{g}.

\begin{lemma} Any nonempty bounded open convex set $K$ in $\mathbf{R}^2$ contains a unique open ellipse $E$ with maximal Euclidean area. Furthermore $\partial \overline{K} \cap \partial \overline{E}$ contains at least three points.
\end{lemma}

Continue the proof of Theorem \ref{1}. We consider the ellipse $E$
in $K$ with maximal Euclidean area. If the number of contact
points of their boundary is 3 or 4. We extend the ellipse $E$ a
little bit to obtain another ellipse $E'$ as in Figure
\ref{extend}.

\begin{figure}[htbp]
\begin{center}
\includegraphics[scale=.4]{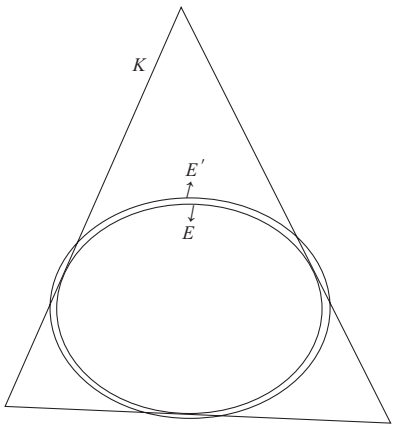}
\end{center}
\caption{\label{extend}Extend ellipse $E$ to $E'$}
\end{figure}

More precisely, assume $E$ is defined by
$$\frac{x_1^2}{a_1^2}+\frac{x_2^2}{a_2^2}<1.$$ We define ellipse $E'$ by
$$\frac{x_1^2}{a_1^2}+\frac{x_2^2}{a_2^2}<1+\varepsilon,$$ where $\varepsilon >0$ and sufficiently small.
Since $E$ is of maximal area, $\partial \overline{E'}$ must intersect
$\partial \overline{K}$ at 6 or 8 points.

Therefore we always can find an ellipse, still denote it by $E$,
such that the number of points in $\partial \overline{K} \cap
\partial \overline{E}$ is at least 5.

Since we assume $K$ is not an ellipse, there exists an open set $U
\subset K-E.$  Now we choose a point $u\in U$ and five points
$p_1,p_2,p_3,p_4,p_5$ from $\partial \overline{K} \cap \partial
\overline{E}$. We arrange the subindex such that the order of
these five points on $\partial \overline{E}$ is
$p_1,p_2,p_5,p_3,p_4$ and $u$ is between $p_1,p_4$ as in Figure
\ref{triangle}. Notice that any three points do not lie on a
straight line.

\begin{figure}[htbp]
\begin{center}
\includegraphics[scale=.4]{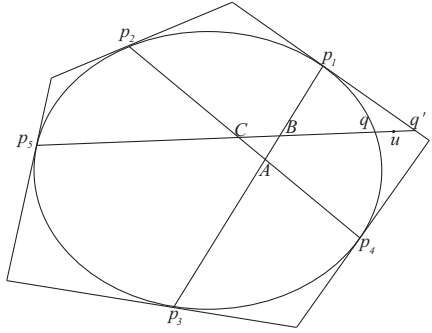}
\end{center}
\caption{\label{triangle}Construction of triangle
$\triangle_{ABC}$}
\end{figure}

Now join $p_1$ with $p_3$, $p_2$ with $p_4$ by line segments,
denoted by $p_1p_3, p_2p_4$. Assume $ p_1p_3\cap p_2p_4=A.$ If the
segment $p_5u$ passes through $A$, we can choose another point
from $U$. Therefore we may assume $p_5u$ avoids $A$ and $p_5u\cap
p_1p_3=B$, $p_5u\cap p_2p_4=C.$

Now we claim that the triangle $\triangle_{ABC}$ in $K$ dose not
have the property M. Since $\triangle_{ABC}$ is also in $(E,d_E)$
which is a hyperbolic geometry. Under the distance $d_E$, three
medians of $\triangle_{ABC}$ all pass through one point. We assume
under the distance $d_E$, points $A',B',C'$ are middle points of
edge $BC,AC,AB$ respectively. By the construction of
$\triangle_{ABC}$, we see under the distance $d_K$ the middle
points of edge $BC,AC,AB$ are $A'',B',C'.$ Now we only need to
check $A'\neq A''$ which will finish the proof. Since
$AA',BB',CC'$ all pass through one point, $AA'',BB',CC'$ can not
all pass through one point.

Now we extend $p_5u$ such that it intersects $\partial
\overline{E}$ at $q$ and $\partial \overline{K}$ at $q'$ as in
Figure \ref{triangle}. Since the cross-ration is invariant under
any projective mapping, consider a projection $P$ sending segment
$p_5u$ to the coordinate line $\mathbf{R}$ such that $P(p_5)=0,
P(C)=1, P(B)=b, P(q)=x, P(q')=x'$ as in Figure \ref{projection}.

\begin{figure}[htbp]
\begin{center}
\includegraphics[scale=.4]{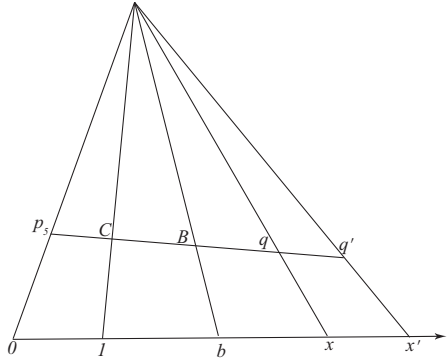}
\end{center}
\caption{\label{projection}Projective mapping}
\end{figure}

Notice $x<x'.$ Furthermore $P(A')=m$ (respectively $P(A'')=m'$)
which is the middle point of $1,b$ with two end points $1,x$
(respectively $1,x'$). In fact, since
$$\frac{m(x-1)}{x-m}=d(1,m)=d(m,b)=\frac{b(x-m)}{m(x-b)},$$ we
have
$$m=\frac{\sqrt{b}x}{\sqrt{b}+\sqrt{(x-1)(x-b)}}=:f(x).$$
We define the right hand side to be a function $f(x).$ And
$m'=f(x').$ Since the derivative $f'(x)<0,$ the function $f(x)$ is
monotonically decreasing in $x.$ Hence $m=f(x)>f(x')=m'$ and
$A'\neq A''.$

For general case $n>2.$ Let us recall the following result, see \cite{bu} (pp. 91).

\begin{lemma}\label{3} For a nonempty open convex set $K$ in $\mathbf{R}^n$, if for a fixed $r, 2\leq r\leq n-1,$
every $r$ dimensional plain through a fixed point $P \in K$
intersects $K$ in an ellipsoid (ellipse when r=2), $K$ itself is
an ellipsoid.
\end{lemma}

Continue the proof of Theorem \ref{1}. If every triangle in $K$
has the property M, by case $n=2$, the intersection of $K$ with
any 2-plain is an ellipse. By the Lemma \ref{3}, $K$ itself is an
ellipsoid.

\end{proof}

\section*{Acknowledgement} The author would like to thank his
advisor, Feng Luo, for suggesting this problem and helpful
discussions.

\bibliographystyle{amsplain}

\end{document}